\documentclass[11pt]{article}
\usepackage{amsmath}
\usepackage{amsfonts}
\usepackage{graphicx}
\usepackage{enumerate}
\usepackage[sort,compress]{cite}
\usepackage[papersize={8.5in,11in},margin=1in]{geometry}
\begin{document}

\begin{center}
\Large\bf{Limitations of polynomial chaos expansions in the Bayesian solution
of inverse problems}
\end{center}
\vspace{.5mm}

\begin{center}
Fei Lu$^{1,2}$, Matthias Morzfeld$^{1,2}$, Xuemin Tu$^3$ and Alexandre J. Chorin$^{1,2}$
\vspace{3mm}

$^1$Lawrence Berkeley National Laboratory;\\
\vspace{1mm}

$^2$Department of Mathematics\\
University of California, Berkeley;\\
\vspace{1mm}

$^3$Department of Mathematics\\
University of Kansas.\\
\end{center}

\begin{abstract}
Polynomial chaos expansions are used to reduce the computational cost in the
Bayesian solutions of inverse problems by creating a surrogate posterior
that can be evaluated inexpensively. We show, by analysis and example, that
when the data contain significant information beyond what is assumed in the
prior, the surrogate posterior can be very different from the 
posterior, and the resulting estimates become inaccurate. One can improve
the accuracy by adaptively increasing the order of the polynomial chaos, 
but the cost may increase too fast for this to be
cost effective compared to Monte Carlo sampling without a surrogate posterior.
\end{abstract}

\section{Introduction}
There are many situations in science and engineering where one needs to estimate parameters
in a model, for example, the permeability of a porous medium in studies of
subsurface flow, on the basis of noisy and/or incomplete data, e.g.~pressure measurements. In the
Bayesian approach, prior information and a likelihood function for the data
are combined to yield a posterior probability density function (pdf) for the
parameters. This posterior can be approximated by Monte Carlo sampling and
in principle yields all the information one needs, in particular the
posterior mean (see e.g.~\cite{KS05,Stu10,Ta05}). 
However, the sampling may require the evaluation of the posterior 
for many values of the parameters, 
which in turn requires repeated solution of the forward problem.
This can be expensive,
especially in complex high-dimensional problems.

Polynomial chaos expansions (PCE) and generalized PCE provide an approximate representation of the solution of the forward problem (see e.g.~\cite{GS03,LK10,Naj09,Xiu10}) which can be used to reduce the cost of Bayesian inverse problems \cite{BRIFG,MNR07,MN09, MX09,
RND12}. The PCE leads to an approximate representation of the posterior,
called a ``surrogate posterior", which can generate a large number of
samples at low cost. However, the resulting samples approximate the surrogate
posterior, not the posterior, so that the accuracy of estimates based on these samples
depends on how well the surrogate posterior approximates the posterior.

We study how the accuracy of the surrogate posterior depends on the data,
and show that when the data are informative (in the sense that the posterior differs
significantly from the prior), then the surrogate posterior can be very different
from the posterior and PCE-based sampling is either inaccurate or
prohibitively expensive. 
Specifically, we examine the behavior of PCE-based sampling  in the small noise regime 
\cite{VW12,VW13}, and report
results from numerical experiments on an elliptic inverse problem for
subsurface flow. In the example, a sufficiently accurate 
PCE require a high order, which makes PCE-based sampling
expensive compared to sampling the posterior directly, 
without a PCE. 
Other limitations of PCE have been 
reported and discussed in other settings as well, e.g. in uncertainty quantification 
\cite{BM13,HLRZ06,BZ14}, and statistical hydrodynamics \cite{Cho74,CC70}. 

The paper is organized as follows. In Section 2 we explain the use of 
PCE in the Bayesian solution of inverse problems. In section 3 we analyze the accuracy of
the surrogate posterior in the small noise regime. In Section 4 we study the efficiency of
PCE-based sampling with numerical examples. 
Section 5 provides a summary.
Proofs and derivations can be found in the appendix.

\section{Polynomial chaos expansion for Bayesian inverse problems}
Consider the problem of estimating model parameters $\theta \in \mathbb{R}%
^{m}$ from noisy data $d\in \mathbb{R}^n$ such that:
\begin{equation}
\label{eq:data}
d=h(\theta )+\eta , 
\end{equation}
where $h:\mathbb{R}^{m}\rightarrow \mathbb{R}^{n}$ is a smooth nonlinear function describing how the parameters affect the data, and where $\eta\sim p_\eta(\cdot)$ is a random variable with  known pdf that represents uncertainty in the measurements. Here, $h$ is the model and often involves a partial differential equation (PDE), or a discretization of a PDE, in which case the evaluation of $h$ can be computationally expensive. 
Following the Bayesian approach, we assume that prior information about the parameters is available in form of a pdf $p_0(\theta)$. This prior and the likelihood $p(d\vert \theta)=p_\eta(d-h(\theta) )$, defined by (\ref{eq:data}), are combined in Bayes' rule to give the posterior pdf
\begin{equation}
p(\theta |d)=\frac{1}{\gamma (d)}p_{0}(\theta )p(d|\theta ),  \label{eq:post}
\end{equation}
where  $\gamma (d)=\int p_{0}(\theta )p(d|\theta )d\theta $ is a normalizing constant (the marginal probability of the data). For simplicity, we assume throughout this paper that  $\eta\sim\mathcal{N}(0,\sigma ^{2}I_n)$ is Gaussian with mean zero and variance $\sigma^2I_n$, and that the prior is $p_0(\theta)=\mathcal{N}(0,I_m)$
(here, $I_k$ is the identity matrix of order $k$). These assumptions may be relaxed,
however we can make our points in this simplified setting. 
In this context, it is important to point out that we make no assumptions about the underlying (numerical) model which, 
in most cases, is nonlinear.

In practice, Monte Carlo (MC) methods such as importance
sampling or Markov chain Monte Carlo (MCMC) are used to represent the posterior numerically (see e.g.~\cite{CH13,Liu01}). Most MC sampling methods require repeated 
 evaluation of the posterior for many instances of $\theta $. Since each posterior
evaluation involves a likelihood evaluation, many evaluations of the model are needed, which can be computationally expensive.

To reduce the computational cost of MC sampling one can approximate the model by a truncation of its PCE, because the evaluation of the truncated PCE is often less expensive than the evaluation of the model (e.g.~solving a PDE). 
It is natural to construct the PCE before the data are available, i.e.~one expands $h$ using the prior. 
With a Gaussian
prior one uses (multivariate) Hermite polynomials, which form a complete
orthonormal basis in $L^{2}(\mathbb{R}^{m},p_{0})$. Let $i=\left(
i_{1},\dots ,i_{m}\right) \in \mathbb{N}^{m}$ be a multi-index and let $%
\theta =(\theta _{1},\theta _{2},\dots ,\theta _{m})$ be the parameter we
wish to estimate. The multivariate Hermite polynomials $\left\{ \Phi _{%
i}(\theta ):\left\vert i\right\vert =i_{1}+\cdots
+i_{m}<\infty \right\} $ are defined by
\begin{equation*}
\Phi _{i}(\theta )=H_{i_{1}}(\theta _{1})\cdots H_{i_{m}}(\theta
_{m}),
\end{equation*}%
where $H_{k}(x)$ is the normalized $k$th-order Hermite polynomial (see e.g.~%
\cite{Sch00,Xiu10}). Assuming that $h\in L^{2}(\mathbb{R}^{m},p_{0})$ we
define the $N$th-order PCE of $h$ by
\begin{equation*}
h_{N}(\theta )=\sum_{\left\vert i\right\vert \leq N}a_{i%
}\Phi _{i}(\theta ),
\end{equation*}%
where the coefficients $a_{i}$ are given by
\begin{equation*}
a_{i}=E[h(\theta )\Phi _{i}(\theta )]=\int h(\theta )\Phi
_{i}(\theta )p_{0}(\theta )d\theta .  \label{ui}
\end{equation*}%
As $N\rightarrow \infty $, $h_{N}$ converges to $h$ in $L^{2}(\mathbb{R}%
^{m},p_{0})$. The rate of convergence depends on the regularity of $h$ and
is estimated by (see e.g. ~\cite{XG03})
\begin{equation}
\Vert h-h_{N}\Vert _{L^{2}(\mathbb{R}^{m},p_{0})}\leq CN^{-\frac{k}{2}}\Vert
h\Vert _{k,2},  \label{convL2}
\end{equation}%
where $C$ is a constant depending only on $m$ and $k$, and $\Vert h\Vert
_{k,2}$ is the weighted Sobolev norm defined by 
$\Vert h\Vert _{k,2}^{2}=\sum_{|\alpha |\leq k}\Vert D^{\alpha }h\Vert
_{L^{2}(\mathbb{R}^{m},p_{0})}^{2}$ 
with $D^{\alpha }h=\frac{\partial ^{|\alpha |}}{\partial _{x_{1}}^{\alpha
_{1}}\cdots \partial _{x_{m}}^{\alpha _{m}}}h$. For the remainder of this paper we assume enough regularity of $h$, so that $\Vert h-h_{N}\Vert _{L^{2}(\mathbb{R}^{m},p_{0})}$ converges quickly with $N$.

In PCE-based sampling for Bayesian inverse problems, one replaces the model $%
h$ in (\ref{eq:data}) by its truncated PCE $h_{N}$, and obtains the
surrogate posterior
\begin{equation}
p_{N}(\theta |d)=\frac{1}{\gamma _{N}(d)}p_{0}(\theta )p_{\eta
}(d-h_{N}(\theta )),  \label{eq:postN}
\end{equation}%
where $\gamma _{N}(d)=\int p_{0}(\theta )p_{\eta }(d-h_{N}(\theta ))d\theta $. 
This surrogate posterior converges to the posterior at the same rate as $\Vert h-h_{N}\Vert _{L^{2}(\mathbb{R}%
^{m},p_{0})}$ converges to zero as $N\rightarrow \infty $ (see equation (\ref{convL2})) 
in the Kullback-Leibler divergence (KLD)
\cite{MX09} and in the Hellinger metric \cite{Stu10}. 

In practice, PCEs of small to moderate order are used because otherwise PCEs
become expensive (see e.g.~\cite{LK10, Xiu10} and section 4).
This truncation introduces error 
unless the problem is linear or
well represented by a low-order polynomial.
If the truncation error is large,
then the surrogate posterior may be 
very different from the posterior. 
The samples one draws with PCE-based sampling methods
approximate the surrogate posterior, 
which implies that the applicability of PCE-based sampling 
for inverse problems depends on how
well the surrogate posterior approximates the posterior, 
i.e.~on the accuracy of the surrogate posterior.

\section{Accuracy of the surrogate posterior}
We wish to study the effects of inaccuracies of a truncated PCE
on PCE-based sampling methods for inverse problems.
Inaccuracies in the PCE are caused by interaction of two mechanisms:
\begin{enumerate}
	\item The error due to truncation is large
	if the physical model $h(\theta)$ in (\ref{eq:data})
	is poorly represented by a low-order polynomial.
	\item The surrogate posterior must be constructed based on the prior (see above).
	However, the posterior can put significant probability mass in parameter regions
	which are unlikely with respect to the prior.
	Thus, if the model is nonlinear, the PCE may be a poor approximation in the region
	of large posterior probability.
	In this case, the surrogate posterior is a poor approximation of the posterior
	in the region where significant posterior probability mass is located.
\end{enumerate}

We assume throughout that $h(\theta)$ is nonlinear
with high-order polynomials in its PCE,
so that a truncated PCE (of moderate order)
is only locally accurate.
What we want to find are the conditions
under which the lack of
global accuracy of the PCE causes PCE-based sampling 
in inverse problems to be inaccurate.

Our analysis tool is the KLD of the posterior (\ref{eq:post}) and the surrogate posterior (\ref{eq:postN}):  
\begin{equation*}
D_{KL}(p_{N}||p):=\int p_{N}(\theta |d)\ln {\frac{p_{N}(\theta |d)}{p(\theta
|d)}}d\theta.
\end{equation*}
Since the posterior and the surrogate posterior depend on the data, $D_{KL}(p_{N}||p)$ also depends on the data. 
We fix the data, as well as the order $N$ of the truncation (assuming non-zero truncation error),
and show that the surrogate posterior is inaccurate if the data are informative, 
i.e.~if the likelihood moves the posterior away from the prior. 
We consider the small noise regime introduced in \cite{VW12,VW13}, 
where the variance of the prior and the likelihood are small. 
The small noise regime is important in data assimilation because it corresponds
to a case where the posterior probability mass localizes in a ``small''
region of the parameter space. 
Moreover, the small noise regime allows for a rigorous analysis and the results
can be indicative of more general situations. 
Specifically, we pick a prior $p_{0,\epsilon}$ with mean zero and variance $\epsilon$ 
and set $\eta\sim\mathcal{N}(0,\sigma^2\epsilon I)$. 
These choices give the small noise posterior
\begin{equation*}
	p_{\epsilon }=\frac{1}{\gamma _{\epsilon }(d)}\exp\left(-\frac{1}{\epsilon }F(\theta )\right),  \label{postEP}
\end{equation*}
where
\begin{equation*}
	F(\theta )=\frac{||d-h(\theta )||^{2}}{2\sigma ^{2}}+\frac{||\theta ||^{2}}{2},
\end{equation*}
and where $\gamma _{\epsilon }(d)=\int \exp\left(-F(\theta )/\epsilon\right)d\theta$ is the normalizing constant; 
here, and for the remainder of this paper $\vert\vert\cdot\vert\vert$ is the Euclidean norm.
We now state our two main results.

\vspace{5mm}
{\bf{Claim 1}}. The KLD of the posterior and the prior grows as $\epsilon$ 
 becomes smaller and smaller. More precisely, we have
\begin{equation}
D_{KL}(p_{0,\epsilon }||p_{\epsilon })=\frac{1}{\epsilon }(F(0)-\min_{\theta
}F(\theta))+o(\frac{1}{\epsilon }),  \label{KLDpri-post}
\end{equation}
which grows to infinity as $\epsilon \rightarrow 0$ when $F(0)>\min_{\theta
}F(\theta)$; here we use the standard notation $f(\epsilon )=o(1/\epsilon)$ 
for $\lim \sup_{\epsilon \rightarrow 0}\epsilon f(\epsilon )=0$. 

\vspace{5mm}
{\bf{Claim 2}}. The KLD of the surrogate posterior to the posterior
\begin{equation}
D_{KL}(p_{N,\epsilon }||p_{\epsilon })=\frac{1}{\epsilon }(F(\theta
_{N}^{\ast })-\min_{\theta }F(\theta ))+o(\frac{1}{\epsilon }),
\label{KLDpc-post}
\end{equation}
grows to infinity as $\epsilon \rightarrow 0$ if $F(\theta _{N}^{\ast
})>\min_{\theta }F(\theta )$; here, $\theta _{N}^{\ast }=\arg
\min_{\theta}F_{N}(\theta )$ with
\begin{equation*}
F_{N}(\theta )=\frac{||d-h_{N}(\theta )||^{2}}{2\sigma ^{2}}+\frac{||\theta
||^{2}}{2}
\end{equation*}
and we assume that the minimizer $\theta _{N}^{\ast }$ is unique.

\vspace{5mm}
Derivations of these two results are provided in the appendix.
The interpretation of the above results is as follows.
Claim 1 shows that, under our assumptions of small noise,
the data become more informative as $\epsilon\to 0$,
because the data shifts the probability mass away from the prior.
Claim 2 thus shows that the surrogate posterior diverges from 
the posterior as the data become more informative.
The two claims combined show that the surrogate posterior 
may not be useful when the data are informative.
We will study the effects of these inaccuracies on the solution 
of Bayesian inverse problems with numerical examples in the next section.

Our results can also be interpreted geometrically. 
As $\epsilon$ is getting smaller,
the posterior is getting more sharply peaked around its mode, since, 
from (\ref{postEP}) and (\ref{Gamma_e}) one obtains
\begin{equation}
p_{\epsilon }(\theta |d)=\exp\left(-\frac{1}{\epsilon }(F(\theta )-\min_\theta F)+o(\frac{1}{\epsilon })\right).\label{postEP2}
\end{equation}
Similarly, one can show that the surrogate posterior is also sharply peaked around its maximum.
Thus, when the maxima of the posterior and the surrogate posterior are different, 
the surrogate posterior is (almost) singular with respect to the posterior.

\section{Efficiency of the PCE for inverse problems}
We have shown that, for a fixed $N$, the KLD of the surrogate posterior from the posterior can be large, i.e.~the PCE is not very accurate. To make it accurate one can increase the order $N$ of the truncated PCE
(in fact, for $N\to\infty$, the PCE is exact everywhere). What we must find out is the rate at which $N$ must be increased to ensure sufficient accuracy. If this rate is too rapid, the PCE becomes increasingly expensive. For example, a  stochastic collocation routinely requires (at least) $N+1$ quadrature points per parameter to compute the coefficients of the PCE, so that the cost of constructing a PCE of order $N$ is about $(N+1)^m$. If  accuracy requires large $N$, then PCE may no longer be a cost effective approach to the inverse problem. 

We can estimate the rate at which $N$ must increase from (\ref{KLDpc-post}), which indicates that the minimizer $\theta _{N}^{\ast }$ of $F_{N}$ must be close to the minimizer $\theta ^{\ast }$ of $F$, or else the KLD, and, hence, the errors are large. Thus, point-wise accuracy of the truncated PCE is needed at least up to $\theta ^{\ast }$, i.e.~making $\sup_{x\in B}|h(x)-h_{N}(x)|$ small for a ball $B\in \mathbb{R}^{m}$ centered at zero and containing $\theta ^{\ast }$. The estimate
\begin{equation*}
\sup_{x\in \mathbb{R}^{m}}\exp\left(-\frac{||x||^{2}}{4}\right)||h(x)-h_{N}(x)||\leq CN^{%
\frac{m}{4}-\frac{k}{2}}\Vert h\Vert _{k,2}. 
\end{equation*}%
from \cite{XG03} indicates that this point-wise accuracy can require that
\begin{equation}
\label{eq:N}
	N>C\exp\left(\frac{||\theta ^{\ast }||^{2}}{(2k-m)}\right).
\end{equation}
Recall that $h$ is smooth, so that $2k>m$ (i.e.~the exponent is positive).
Moreover, since the mean of the prior is zero, a large $\theta ^\ast$ (in Euclidean norm) is far from where the prior probability mass is. Thus, large $\theta^*$ indicates that the data are informative. Equation~(\ref{eq:N}) thus shows that the order required to maintain accuracy in the PCE grows quickly as the data become increasingly informative.

We now investigate the effects of the inaccuracy of the surrogate posterior on the efficiency of PCE-based sampling for inverse problems, using numerical experiments with an elliptic inverse problem. We choose an elliptic inverse problem because it is a popular tool for testing algorithms for parameter estimation and it is also theoretically well understood \cite{MN09,MX09,MNR07,TMWC,DS11,Stu10}. The example is not very realistic, however it is sufficient to help us make our points.

\subsection{Numerical experiments}
Consider the random elliptic equation
\begin{equation}
\left\{
\begin{array}{ll}
-\nabla .(e^{Y(x)}\nabla u)=f(x), & x\in D; \\
u(x)=0, & x\in \partial D,%
\end{array}%
\right.   \label{epde}
\end{equation}%
where $x=(x_1,x_2)$, $D=\left[0,1\right]\times\left[0,1\right]$, $f(x)=\pi ^{2}\sin (\pi x_{1})\sin (\pi x_{2})$, and where $Y$ (often called the log-permeability) is a square integrable stochastic field we wish to estimate from data (e.g.~noisy measurements of $u$ at some locations in $D$). Equation (\ref{epde}) is a simplified model for flow in porous media \cite{DS11,TMWC}. 
We consider a Gaussian log-permeability  with mean zero and squared exponential covariance function (see, e.g.~\cite{Rasmussen2006})
\begin{equation}\label{equation:correlation}
R(x_1,x_2,y_1,y_2)=\sigma_1^2\sigma_2^2\exp \left(-\frac{\left(x_1-y_1\right)^2}{l_1}-\frac{\left(x_2-y_2\right)^2}{l_2}\right).
\end{equation}
In our numerical experiments we set $\sigma_1=\sigma_2=1$ and the correlation lengths $l_1=l_2=1$.

We discretize (\ref{epde}) by a standard finite element method 
on a uniform $(M+1)\times (M+1)$ mesh of triangular 2-D elements \cite{Braess}, with $M=15$.
Solving the PDE is thus equivalent to solving the linear system
\begin{equation*}
A(\theta )\hat{u}(\theta )=\hat{f},  \label{depde}
\end{equation*}
where $A(\theta)$ is a $M^2\times M^2$ symmetric positive
definite matrix, $\hat{u}(\theta)$ is a $M^2$-dimensional vector
approximating the continuous solution $u$, 
and where $\hat{f}$ is a discrete version of $f$. 
We evaluate integrals by Gaussian quadrature.
We solve these equations with a preconditioned conjugate gradient method. 
The infinite dimensional random field $Y$ is discretized by
a Karhunen--Lo\`{e}ve expansion (see, e.g.~\cite{GS03,LK10})
\begin{equation*}
	\hat{Y}=U\theta,
\end{equation*} 
where $\theta =(\theta _{1},\dots ,\theta _{m})$ and where $U=(U_{1}(x),\dots ,U_{m}(x))$ is a matrix whose columns are the first $m$ eigenvectors of the squared exponential covariance function, multiplied by their corresponding eigenvalue. Since the squared exponential covariance function has a rapidly decaying spectrum, we can capture $96.66\%$ of the variance with $m=3$, and this is the choice we make for the following numerical experiments.

This set-up implies a Gaussian prior for $\theta$, with mean zero and variance $I$. 
We expand the solution $\hat{u}(\theta )$ of (\ref{depde}) in a PCE and approximate it by the truncated PCE, $\hat{u}_{N}(\theta )$, of order $N$. We use stochastic collocation with $N+2$ ($N$ even) Gaussian--Hermite quadrature points per parameter to compute the coefficients of the PCE. Thus, $(N+2)^m$ PDE solutions are required to construct the PCE. This approach is only efficient if $N$ or $m$ are small. Here, $m=3$ is small, and we study the effects of the order $N$ of the truncation, i.e.~we will vary $N$ in our numerical experiments.

The example therefore corresponds to a situation in which PCE could be useful,
because the number of dimensions is small (it is equal to 3). 
If we decrease the correlation length, which is perhaps more realistic,
we would need to increase $m$ to capture the variance
and PCE becomes impractical.
However, our goal is to show how inaccuracies in the surrogate posterior
can be caused by a PCE which is locally a good approximation,
but which lacks global accuracy.

\begin{figure}[tb]
\centering\includegraphics[width=0.8\textwidth]{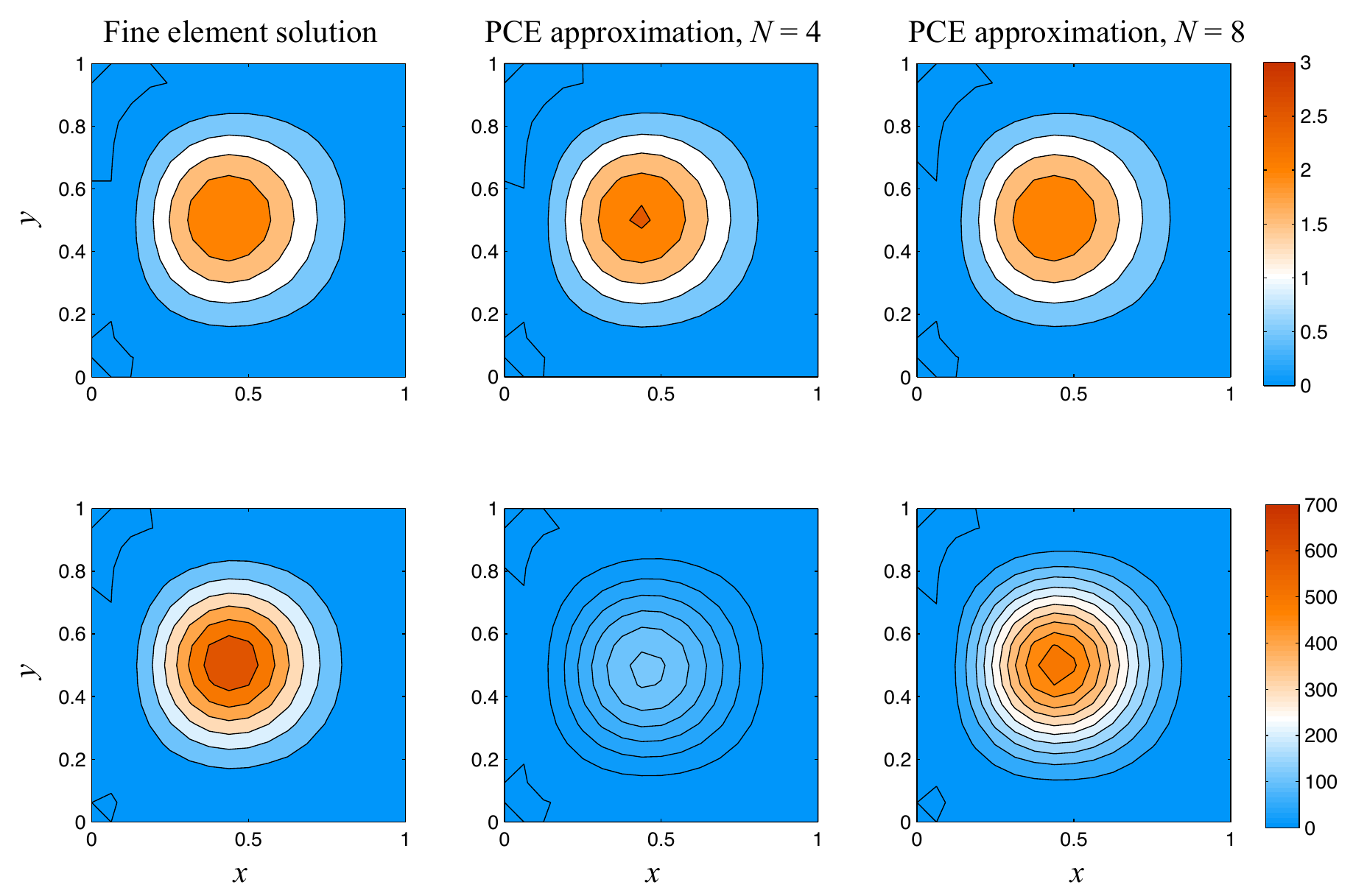}
\caption{ The finite element solutions (left column), and their PCE
approximations of order $N=4$ (middle column) and of order $N=8$ (right
column) for two different parameters: $\theta =(-2,-1,1)$ (top row) and $\protect\theta %
=(-8,-1,1)$ (bottom row).}
\label{fig_solu_contour}
\end{figure}

Figure~\ref{fig_solu_contour} shows our finite element solutions and their PCE
approximations of order $N=4$ and $N=8$, for two different values of the parameter $\theta$.
The first parameter, $\theta =(-2,-1,1)$, is close to the mean of the prior and we observe that both PCE approximations are accurate (see the top row of Figure~\ref{fig_solu_contour}). As we move the parameter further away from the mean of the prior, e.g.~choose $\theta =(-8,-1,1)$, we observe that the accuracy of the PCE requires that $N\geq8$  (see the bottom row of Figure \ref{fig_solu_contour}).

We study the accuracy of the PCE approximation further by focusing attention on the grid point $x=(0.5,0.5625)$,
which is the point on our grid where the first eigenvector of the squared exponential covariance function is maximum. 
We vary $\theta_1$ between $-8,\dots,2$ and fix $\theta_2=-1$ and $\theta_3=1$ as before, i.e. one parameter, $\theta_1$, may be far from its value assumed by the prior, while the other two are within two standard deviations of the mean of the prior. We restrict $\theta _{1}\leq 2$ because the finite element solution is very small otherwise (less than the  standard deviation of the noise in the data), which corresponds to a situation in which no inference can be successful since there is almost no information in the signal. The results of our calculations are shown in table 1, where we show the relative error as a function of $\theta_1$ for PCEs with $N=4$ and $N=8$. 
\begin{table}[tb]
\begin{center}
\label{tab:Results}
\begin{tabular}{ccccccccccccc}
 Order     & \multicolumn{11}{c}{Parameter $\theta_1$}\\
      
& & -8 & -7 & -6 & -5 & -4 & -3 & -2 & -1 & 0 & 1 & 2\\
\hline
 $N=4$ & 81.1&    69.4&    53.2&  33.4&    13.5&   0.26&   2.56&  2.68&   1.71&   7.72&    85.3\\
 $N=8$ & 21.2&    10.6&    3.51&  0.39&    0.13&   0.03&   0.02&  0.01&   0.02&   0.04&    0.02\\
 \hline
\end{tabular}
\end{center}
\caption{Relative error in percent as a function of of $\theta_1$ at the grid point $x=(0.5,0.5625)$.}
\end{table}
The relative error is defined by the absolute value of the difference of the finite element solution and the PCE approximation, divided by the absolute value of the finite element solution (all evaluated at $x=(0.5,0.5625)$).

We find that the PCE approximation of order $N=4$ is accurate (in the sense that the error less than 10\%) if $\theta_1$ is within a standard deviation of the mean of the prior (i.e., roughly between -2 and 1). One can extend the region where the PCE is accurate by increasing its order, e.g.~a PCE with $N=8$ is accurate (with errors less than 10\%) for $-6\leq\theta_1\leq 2$. However, for finite $m$, the PCE is always inaccurate for large enough $\theta_1$, as is indicated by large errors in Table~1
(global accuracy can only occur if $h(\theta)$ is a low-order polynomial).

With the PCE in place, we solve the inverse problem of estimating the log-permeability~$\theta$ from noisy measurements of the solution $u$ at several locations in the domain $D$. The data are ``synthetic'', i.e.~we generate data using our numerical model. This has the advantage that the sampling algorithms operate under ideal conditions (since the model is compatible
with the data). We collect data in $n=16$ locations in the center of the domain, i.e.
\begin{equation*}
d=G\hat{u}(\theta )+\eta ,  \label{ddata}
\end{equation*}%
where $G$ is an $n\times M^{2}$ matrix which projects the finite element solution to the $n$ observed components, and where $\eta $ is an $n$--dimensional random vector with distribution $\mathcal{N}(0,0.05^{2}I)$. 
More precisely, data are collected
3 steps away from the top and right boundary and five steps away from the bottom and left boundaries, 
and two steps away from each other.
The prior and likelihood define the posterior
\begin{equation*}
p(\theta |d)\propto \exp\left(-F(\theta )\right),  \label{posterior}
\end{equation*}%
where 
\begin{equation}
F(\theta )=\frac{1}{2}||\theta ||^{2}+\frac{1}{2\sigma ^{2}}(d-G\hat{u%
}(\theta ))^{T}(d-G\hat{u}(\theta )). \label{eq:F}
\end{equation}
For PCE-based MC sampling we replace $h$ with its truncated PCE (we consider $N=4$ and $N=8$) and compute the surrogate posterior
\begin{equation*}
p_{N}(\theta | d) \propto \exp\left(-F_N(\theta)\right),
\end{equation*}
where 
\begin{equation*}
	F_{N}(\theta )=\frac{1}{2}||\theta|| ^{2}+\frac{1}{2\sigma ^{2}}(d-G\hat{u}_{N}(\theta ))^{T}(d-G\hat{u}_{N}(\theta )).
\end{equation*}

We generate synthetic data sets in which we vary $\theta_1=-10,-9,\dots ,2$, while $\theta_2=-1$, $\theta_3=1$ are kept fixed so that the data become increasingly informative as $\theta_1$ increases in magnitude and use implicit sampling \cite{AMC13,CMT10,CT09,MTAC12} for each data set to sample the surrogate posterior. Implicit sampling is an importance sampling
method that generates samples close to the mode of the posterior by first locating the mode  via numerical optimization, and then solving data dependent algebraic equations to generate samples in the vicinity of the mode. The weights of the samples are then computed from the importance density which is implicitly defined by these algebraic equations. For further details, also about the implementation of implicit sampling, see \cite{AMC13,CMT10,CT09,MTAC12}. 

\begin{figure}[tb]
\centering\includegraphics[width=0.4\textwidth]{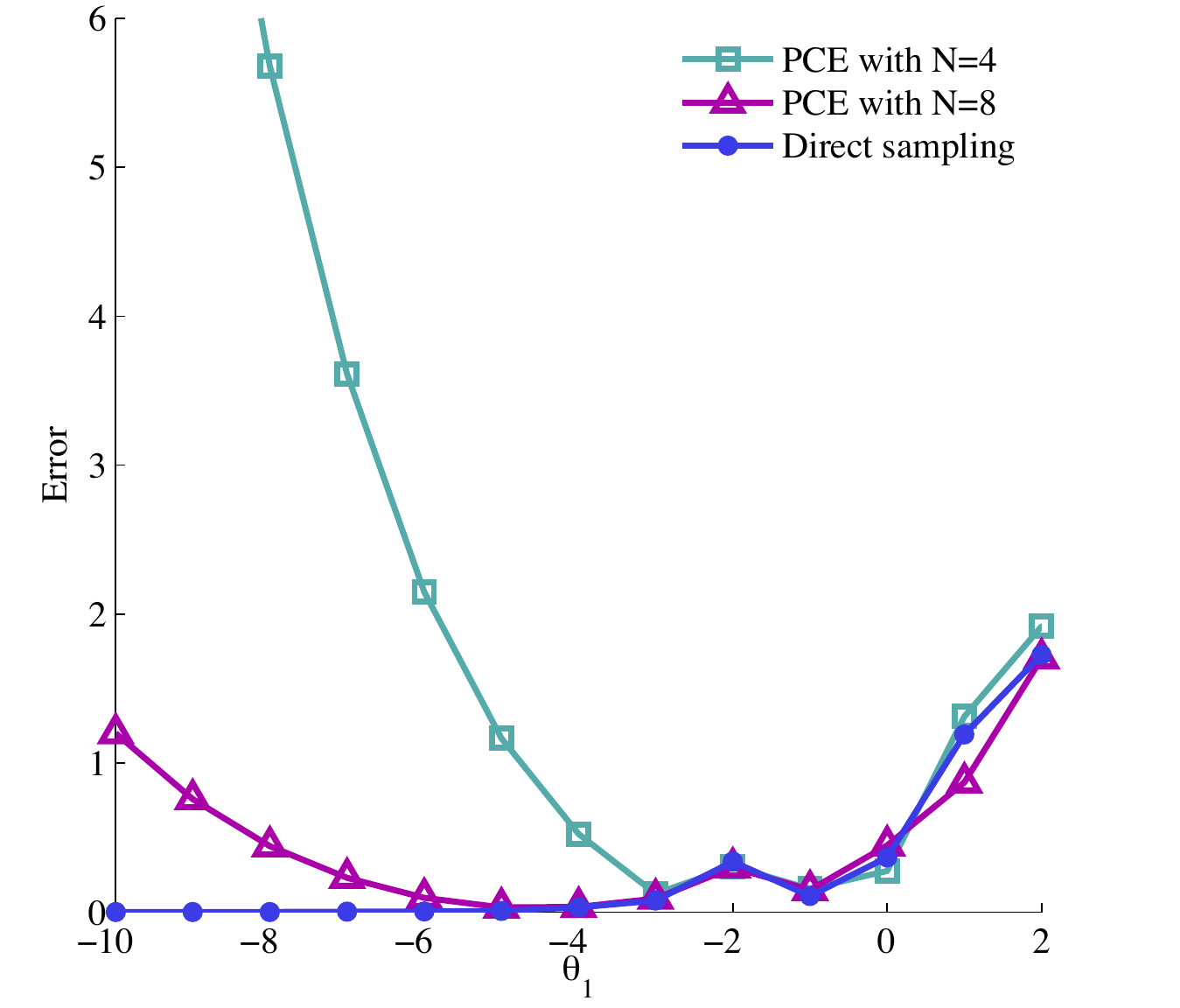}
\caption{Errors in the solution of the inverse problem: norm of the error as a function of $\theta_1$ for implicit sampling (blue dots), PCE-based implicit sampling with $N=4$ (turquoise squares), and PCE-based implicit sampling with $N=8$ (purple triangles). }
\label{fig_err2DM3}
\end{figure}

Figure~\ref{fig_err2DM3} summarizes the results of our numerical experiments.
We plot the error one makes when solving the inverse problem for three different sampling schemes. The error is defined as the Euclidean norm of the difference between the parameter $\theta$ and its approximation by the conditional mean (the minimum mean square error estimator) obtained via sampling. We use implicit sampling with 20 particles, and vary the accuracy of the PCE from $N=4$, to $N=8$, to $N= \infty$, i.e.~we sample the posterior directly (without using PCE). The latter is the reference solution.

The figure illustrates that PCE-based sampling with $N=4$ can only be accurate where the PCE itself is accurate, i.e.~only if the parameter is within 2 standard deviations of the mean of the prior. The error increases steeply as the parameter becomes larger (in magnitude). The figure also indicates that the region of applicability of PCE-based sampling can be increased by increasing the order of the PCE, since we obtain much better results with a surrogate posterior when $N=8$. In this case, the parameter can be far from what is assumed by the prior, however ultimately, one can not guarantee the accuracy of this approach due to the lack of global accuracy of the PCE.

We have obtained the same results with a random walk Metropolis MCMC algorithm (see e.g.~\cite{Liu01}),
which also shows that the failure of PCE we observe here is independent of the sampling method one uses 
for sampling the (surrogate) posterior.

\begin{figure}[tb]
\centering\includegraphics[width=\textwidth]{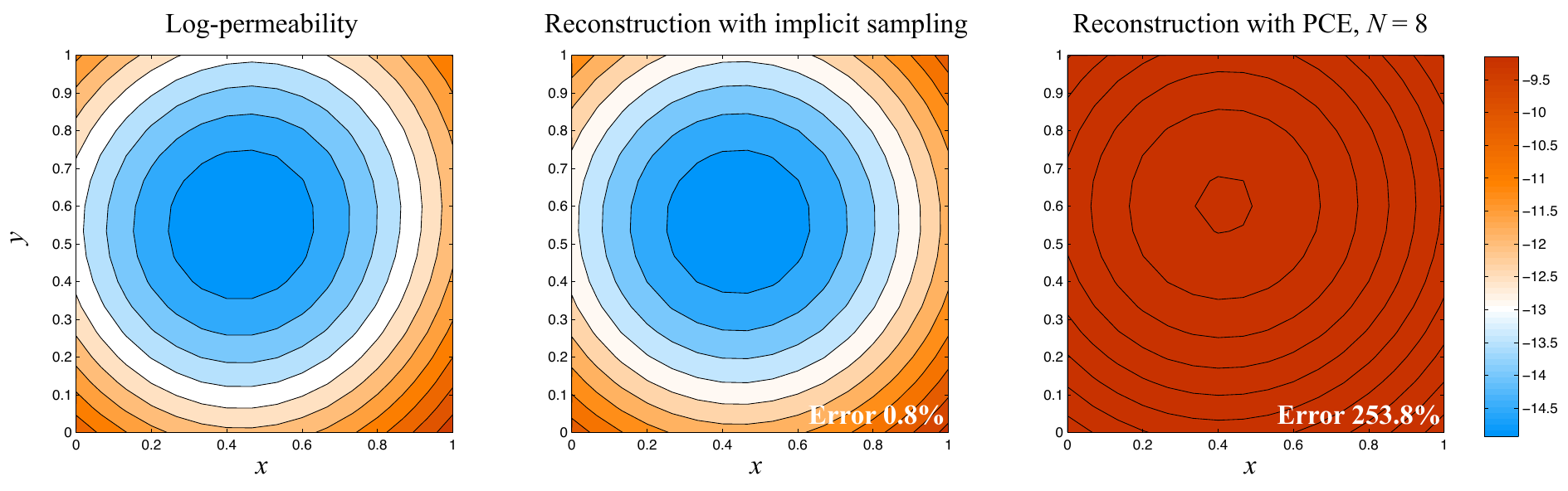} \
\caption{The log-permeability (left), its approximation via implicit sampling (middle) and its approximation via PCE-based MCMC sampling with order $N=8$.}
\label{fig_perm}
\end{figure}

In Figure~\ref{fig_perm} we indicate what the errors in Figure~\ref{fig_err2DM3} mean for the physics of this inverse problem. 
In the left panel we plot the log-permeability, in the middle panel its approximation obtained via implicit sampling of the posterior ($N= \infty$), and in the right panel we plot the estimated log-permeability obtained with PCE-based sampling $N=8$. The parameter $\theta_1=-15$, i.e.~far from what is assumed in the prior, while $\theta_2=-1$ and $\theta_3=1$ are within the range predicted by the prior. This is a scenario in which the data are informative and move the posterior away from the prior. It is evident from this figure that using the surrogate posterior gives catastrophically large errors. However, it is feasible to solve this inverse problem by sampling the posterior 
(without constructing a surrogate posterior).

Finally, we compare the computational costs of MC sampling with and without a PCE-based surrogate posterior. 
The cost of solving the inverse problem is dominated by the cost of the required PDE solutions.
Constructing a PCE with stochastic collocation with $N+2$ Gaussian--Hermite quadrature points for each parameter dimension, requires (at least) $(N+2)^m$ PDE solutions. Since $m=3$ and $N=4$ or $N=8$, between 216 and 1000 PDE solutions are required for constructing the PCE.
 Implicit sampling of the posterior requires 153 evaluations of the posterior, which amounts to 153 PDE solutions (most of them occur during the optimization to find the mode of the posterior, since 20 samples appear sufficient to accurately describe it). 
Neglecting the cost of evaluation of the surrogate posterior, this is a larger cost than implicit sampling of the posterior, however the results we obtain are less accurate. Since PCE-based sampling is more costly, but less accurate, we conclude that sampling the posterior (without constructing a surrogate posterior) is a better method for this~problem.

\subsection{Discussion}
In the small noise regime, the approximation of the posterior in (\ref{postEP2}) indicates that significant information resides in the neighborhood of the maximum of the posterior. Hence a successful sampling method should generate samples around this maximum, otherwise information will be lost. Implicit sampling is such a sampling method and therefore is efficient in the above example
which corresponds to a small noise scenario. In other situations, other sampling schemes may show a better performance in terms of balancing computational cost and accuracy.

Our simulations also suggest that a successful sampling scheme could result from  an adaptive construction of a PCE so that the surrogate
posterior is close to the posterior near the maximizer of the posterior. For instance, one could compute a PCE with respect to $\mathcal{N}(\mu,H^{-1})$, where $\mu $ is the maximizer of the posterior and $H$ is the Hessian of $F$ in (\ref{eq:F}) at the maximum.
This construction can gain efficiency and accuracy over the prior-based surrogate, since PCEs with low to moderate order may be locally sufficiently accurate.
However, such a scheme adds the cost of the optimization to the cost of 
the PCE construction (neglecting the cost of using the PCE during sampling). 
It is not clear to us that this strategy will be more efficient than sampling the posterior directly, since e.g.~implicit sampling can generate samples close to the mode of the posterior at a low cost once the mode is located \cite{AMC13,CMT10,CT09,MTAC12}. 

Last, we wish to point out that one can anticipate similar problems and modes of failure with other reduced order modeling techniques for sampling and inverse problems, since the failures we describe are due to the lack of global accuracy, which is common to all reduced order models.

\section{Conclusions}
We have suggested possible mechanisms of failure of PCE-based
sampling in Bayesian inverse problems. In particular, we showed that if the data contain information beyond what is assumed in the prior, then PCE-based sampling can lead to catastrophically large errors or require excessive computational cost. 
The reason is that PCEs of finite order are not globally accurate (unless the model itself is a low-order polynomial). We presented a rigorous analysis of the failure in the small noise limit, which is characterized by a prior and a likelihood that have ``small'' variances. We also investigated the efficiency of PCE-based sampling in numerical experiments with an elliptic inverse problem. We observed that a sufficiently accurate PCE for this problem requires a high order, which makes the approach impractical compared to directly sampling the posterior (without constructing a PCE). Moreover, even at a low accuracy, PCE-based sampling was found to be more costly than sampling the posterior without a PCE.

\section*{Acknowledgements}
This work was supported in part by the Director, Office of Science,
Computational and Technology Research, U.S. Department of Energy under
Contract No.~DE-AC02-05CH11231, 
and by the National Science Foundation under
grants DMS-1217065, DMS-0955078, and DMS-1115759.

\section*{Appendix}
\subsection*{Derivation of Claim 1}
To prove (\ref{KLDpri-post}), note that
\begin{eqnarray}
D_{KL}(p_{0,\epsilon }||p_{\epsilon }) &=&\int \left( 2\pi \epsilon \right)
^{-\frac{m}{2}}\exp\left(-\frac{||\theta|| ^{2}}{2\epsilon }\right)\ln \left( \frac{\gamma
_{\epsilon }(d)}{\left( 2\pi \epsilon \right) ^{\frac{m}{2}}}\exp\left(\frac{%
||d-h(\theta )||^{2}}{2\epsilon \sigma ^{2}}\right)\right) d\theta,  \notag \\
&=&\ln \frac{\gamma _{\epsilon }(d)}{\left( 2\pi \epsilon \right) ^{\frac{m}{%
2}}}+\frac{1}{\epsilon }E\left[ \frac{||d-h(X_{\epsilon })||^{2}}{2\sigma ^{2}}%
\right] ,  \label{KLDpri-post2}
\end{eqnarray}
where $X_{\epsilon }\sim\mathcal{N(}0,\epsilon I_{m})$. The normalizing
constant $\gamma _{\epsilon }(d)$ can be written as
\begin{align}
\gamma _{\epsilon }(d)=&\,\int \exp\left(-\frac{1}{\epsilon }F(\theta )\right)d\theta =
\left( 2\pi \epsilon \right) ^{\frac{m}{2}}E\left[ \exp\left(-\frac{1}{\epsilon }%
\frac{||d-h(X_{\epsilon })||^{2}}{2\sigma ^{2}}\right)\right] .  \label{Gamma_e}
\end{align}%
Laplace's principle (see e.g.~\cite{Var84}) indicates that
\begin{equation*}
\lim_{\epsilon \rightarrow 0}\ln E\left[ \exp\left(-\frac{1}{\epsilon }\frac{%
||d-h(X_{\epsilon })||^{2}}{2\sigma ^{2}}\right)\right] =\min_\theta F(\theta),
\end{equation*}
which can be written as
\begin{equation*}
E\left[ \exp\left(-\frac{1}{\epsilon }\frac{||d-h(X_{\epsilon })||^{2}}{2\sigma ^{2}}\right)%
\right] =\exp\left(-\frac{1}{\epsilon }\min_\theta F(\theta) +o(\frac{1}{\epsilon }%
)\right).  \label{LDP}
\end{equation*}%
Substituting this equality into (\ref{Gamma_e}), we have
\begin{equation}
\gamma _{\epsilon}(d)=\left(2\pi \epsilon \right) ^{\frac{m}{2}}\exp\left(-\frac{1}{\epsilon }
\min_\theta F(\theta) +o(\frac{1}{\epsilon })\right).  \label{Gamma_e2}
\end{equation}%
Since 
\begin{equation*}
  E\left[\frac{||d-h(X_{\epsilon })||^{2}}{2\sigma
^{2}}\right]\rightarrow \frac{||d-h(0)||^{2}}{2\sigma ^{2}}=F(0)
\end{equation*}
as $\epsilon \rightarrow 0$, we can write the second term in (\ref{KLDpri-post2}) as
\begin{equation*}
\frac{1}{\epsilon }E\left[ \frac{||d-h(X_{\epsilon })||^{2}}{2\sigma ^{2}}%
\right]= \frac{1}{\epsilon } F(0)+ o(\frac{1}{\epsilon }).
\end{equation*}%
Then (\ref{KLDpri-post}) follows by substituting the above equality and (\ref{Gamma_e2}) into (\ref{KLDpri-post2}).

\subsection*{Derivation of Claim 2}
To prove (\ref{KLDpc-post}), express the surrogate posterior as
\begin{equation*}
p_{N,\epsilon }=\frac{1}{\gamma _{N,\epsilon }(d)}\exp\left(-\frac{1}{\epsilon }%
F_{N}(\theta )\right),
\end{equation*}
where $\gamma _{N,\epsilon }(d):=\int \exp\left(-F_{N}(\theta
)/\epsilon\right)d\theta $. The definition of the KLD then gives
\begin{eqnarray}
D_{KL}(p_{N,\epsilon }||p_{\epsilon }) &=&\int 
p_{N,\epsilon}(\theta)
\ln \left( \frac{\gamma
_{\epsilon }(d)}{\gamma _{N,\epsilon }(d)}\exp\left(\frac{1}{\epsilon }\left(
F(\theta )-F_{N}(\theta )\right) \right)\right) d\theta,  \notag \\
&=&\ln \frac{\gamma _{\epsilon }(d)}{\gamma _{N,\epsilon }(d)}+\frac{1}{%
\epsilon }\int p_{N,\epsilon }(\theta )\left( F(\theta )-F_{N}(\theta
)\right) d\theta .  \label{KLDpc2}
\end{eqnarray}%
As before (see (\ref{Gamma_e2})), we have that
\begin{equation}
\gamma_{N,\epsilon }(d)=\left(2\pi \epsilon \right) ^{\frac{m}{2}}\exp\left(-\frac{1}{%
\epsilon }\min_{\theta}F_{N}(\theta)+o(\frac{1}{\epsilon })\right). \label{Gamma_e3}
\end{equation}
Thus,
\begin{equation*}
p_{N,\epsilon }= \exp\left(-\frac{1}{\epsilon }(F_{N}(\theta )-\min F_{N})+o(\frac{1%
}{\epsilon })\right),
\end{equation*}
converges to the delta function $\delta _{\theta _{N}^{\ast }}(\theta
) $ as $\epsilon \rightarrow 0$. It follows that
\begin{equation*}
\lim_{\epsilon \rightarrow 0}\int p_{N,\epsilon }(\theta )\left( F(\theta
)-F_{N}(\theta )\right) d\theta =F(\theta _{N}^{\ast
})-\min_{\theta}F_{N}(\theta),
\end{equation*}
which implies that the second term in (\ref{KLDpc2}) can be written as 
$(F(\theta _{N}^{\ast})-\min_{\theta}F_{N}(\theta))/\epsilon+o(1/\epsilon)$.
Then equation (\ref{KLDpc-post}) follows by substituting (\ref{Gamma_e2}) and (\ref{Gamma_e3}) into (\ref{KLDpc2}).

\bibliographystyle{plain}
\bibliography{references}

\end{document}